\documentclass[journal]{IEEEtran}
\ifCLASSINFOpdf
   \usepackage[pdftex]{graphicx}
\else
   \usepackage[dvips]{graphicx}
\fi
\ifCLASSOPTIONcompsoc
  \usepackage[caption=false,font=normalsize,labelfont=sf,textfont=sf]{subfig}
\else
  \usepackage[caption=false,font=footnotesize]{subfig}
\fi
\usepackage[T1]{fontenc} % optional
\usepackage{amsmath}
\usepackage[cmintegrals]{newtxmath}
\usepackage{bm} % optiona
\usepackage{subfig}
\usepackage{graphics}

\begin{document}
\title{Comments on ``Fractional Extreme Value Adaptive Training Method: Fractional Steepest Descent Approach''}

\author{Abdul~Wahab and Shujaat Khan% <-this % stops a space
%\thanks{Corresponding Author: A. Wahab, Tel.: +82-42-350-4360, Fax: +82-42-350-4310.}%
\thanks{A. Wahab is with NUTECH School of Applied Sciences and Humanities, National University of Technology, Sector I-12, Main IJP Road, 44000, Islamabad, Pakistan (e-mail:  wahab@nutech.edu.pk).}% <-this % stops a space
\thanks{S. Khan is with Bio-imaging, Signal Processing and Learning (BISPL) Laboratory, Department of Bio and Brain Engineering, Korea Advanced Institute of Science and Technology, 34141, Daejeon, South Korea (e-mail: shujaat@kaist.ac.kr).}
\thanks{This work was supported by the Korea Research Fellowship Program through the National Research Foundation (NRF) funded by the Ministry of Science and ICT (NRF- 2015H1D3A1062400).}%
\thanks{Manuscript received December 15, 2017.}}


\markboth{IEEE Transactions on Neural Networks and Learning Systems,~Vol.~-, No.~-, December~2017}%
{Wahab and Khan: Comments on ``Fractional Extreme Value Adaptive Training Method: Fractional Steepest Descent Approach''}

\maketitle

\begin{abstract}
In this comment, we raise serious concerns over the derivation of the rate of convergence of fractional steepest descent algorithm in Fractional Adaptive Learning (FAL) approach presented in ``Fractional Extreme Value Adaptive Training Method: Fractional Steepest Descent Approach'' [IEEE Trans. Neural Netw. Learn. Syst., vol. 26, no. 4, pp. 653--662, April 2015]. We substantiate that the estimate of the rate of convergence is grandiloquent. We also draw attention towards a critical flaw in the design of the algorithm stymieing its applicability for broad adaptive learning problems. Our claims are based on analytical reasoning supported by experimental results.
\end{abstract}

\begin{IEEEkeywords}
Fractional calculus, fractional differential, fractional energy norm, fractional extreme point, fractional gradient.
\end{IEEEkeywords}

% \fi
\IEEEpeerreviewmaketitle

\section{Introduction}

\IEEEPARstart{T}{he} Least Mean Squares (LMS) algorithm is a widely used tool in adaptive signal processing due to its stable performance and simple implementation. However, its convergence is slow. Accordingly, many variants of LMS have been proposed in recent years in order to achieve an accelerated convergence without compromising on the steady-state residual error.  In the same spirit, the FAL method based on a \emph{fractional steepest descent approach} was proposed in \cite{paper}. Unfortunately, the rate of convergence of the FAL algorithm is derived in terms of an approximation of the general update rule that furnishes unreliable estimate. We elaborate on this issue in Section \ref{ss:error}.  Further, we draw attention towards a critical flaw in the design of the algorithm stymieing its applicability on general adaptive learning problems in Section \ref{ss:limitation}.  The consequences of these flaws on the proposed method are discussed in Section \ref{s:consequences}. A brief conclusion is provided in Section \ref{s:conclusion}.

\section{Main Remarks}\label{s:mistake}
In order to facilitate ensuing discussion, we follow the notation and equation numbering used in \cite{paper}, the corrected and  the new numbers are distinguished by a superposed asterisk and a prime, respectively. 

\subsection{Remarks on Convergence Analysis}\label{ss:error}

In \cite{paper}, the update equation of the proposed FAL algorithm based on fractional gradient descent is provided in \eqref{eq:19} as 
\begin{align}
\label{eq:19}
s_{k+1}=s_k-\frac{2\mu\eta}{\Gamma(3-\nu)}\left(s_k-s^{\nu *}\right)^2s_k^{-\nu}, \quad\text{if } \nu\neq 1,2,3.
\tag{19}
\end{align}
Since, \eqref{eq:19} is nonlinear, it is intriguing to derive an explicit expression for $s_k$. Towards this end,  $s_k$ is regarded as a discrete sample of a continuous function $s(t)$ at  $t=k$ in \cite{paper}, and \eqref{eq:19} is converted to an ordinary differential equation  (ODE),
\begin{align}
\label{eq:20}
D^1_ts(t)\cong\frac{-2\mu\eta}{\Gamma(3-\nu)\nu (s^{\nu *})^{\nu-1} s(t)}\left[s(t)-s^{\nu *}\right]^2, %\quad\text{if } \nu\neq 1,2,3,
\tag{20}
\end{align}
using a power series expansion of $s^\nu$ about $s-s^{\nu *}$ (furnishing $s^\nu\cong\nu(s^{\nu*})^{\nu-1} s$). Here, $D^1_t$ is the derivative with respect to $t$. The ODE \eqref{eq:20} is solved in \cite{paper} for $s(t)$, thereby furnishing  
\begin{align}
\label{eq:21}
s_k\cong s^{\nu *} + e^{\left(\frac{-2\mu\eta k}{\Gamma(3-\nu)\nu\left(s^{\nu *}\right)^{\nu-1}}\right)}, \quad\text{if } \nu\neq 1,2,3.
\tag{21}
\end{align}

We argue that the expression \eqref{eq:21}, on which the entire convergence analysis is based, is an \emph{unreliable approximation} of the solution to \eqref{eq:20}.  In fact, by separation of variables, \eqref{eq:20} renders
%\begin{align}
%\frac{s(t)}{\left(s(t)-s^{\nu *}\right)^2}D^1_ts(t)\cong \frac{-2\mu\eta}{\Gamma(3-\nu)\nu\left(s^{\nu *}\right)^{\nu-1}},
%\label{eq:a}
%\tag{1'}
%\end{align}
%which on simple integration yields  
\begin{align}
\ln\left |s(t)-s^{\nu *}\right|-\frac{s^{\nu *}}{s(t)-s^{\nu *}}\cong \frac{-2\mu\eta t}{\Gamma(3-\nu)\nu\left(s^{\nu *}\right)^{\nu-1}}+C,
\label{eq:b}
\tag{1'}
\end{align}
where $C$ is the constant of integration whose
 %Here, \eqref{eq:b} is obtained by noting that the right hand side (RHS) of \eqref{eq:a} is constant and for the left hand side (LHS), we have the formula
%\begin{align}
%\int \frac{x}{\left(x-a\right)^2}dx 
%=& \int \frac{1}{\left(x-a\right)}dx +\int \frac{a}{\left(x-a\right)^2}dx
%\nonumber
%\\
%=& \ln |x-a|- \frac{a}{(x-a)}+\widetilde{C},
%\label{eq:c}
%\tag{3'}
%\end{align}
%with any $a\in\mathbb{R}$ and integration constant $\widetilde{C}$. 
value can be determined by the initial input $s_0=s(0)$. Specifically, 
\begin{align}
\label{eq:d}
C\cong \ln\left |s_0-s^{\nu *}\right|-[{s^{\nu *}}/({s_0-s^{\nu *}})].
\tag{2'}
\end{align}
Substituting \eqref{eq:d} in \eqref{eq:b} and setting $s(t)=s_k$, one gets 
%\begin{align*}
%\ln\left |\frac{s(t)-s^{\nu *}}{s_0-s^{\nu *}}\right|\cong \frac{-2\mu\eta t}{\Gamma(3-\nu)\nu\left(s^{\nu *}\right)^{\nu-1}}+\frac{s^{\nu *}}{s(t)-s^{\nu *}}-\frac{s^{\nu *}}{s_0-s^{\nu *}},
%\end{align*}
%or equivalently 
%\begin{align}
%(s(t)-s^{\nu *})\cong\left(s_0-s^{\nu *}\right)\exp\Bigg(&\frac{-2\mu\eta t}{\Gamma(3-\nu)\nu\left(s^{\nu *}\right)^{\nu-1}}
%\nonumber
%\\
%&+\frac{s^{\nu *}}{s(t)-s^{\nu *}}-\frac{s^{\nu *}}{s_0-s^{\nu *}}\Bigg).
%\label{eq:e}
%\tag{5'}
%\end{align}
%Based on this expression, one finds out that
\begin{align}
(s_k-s^{\nu *})
&\cong\left(s_0-s^{\nu *}\right)e^{\left(\frac{-2\mu\eta k}{\Gamma(3-\nu)\nu\left(s^{\nu *}\right)^{\nu-1}}\right)}
%\nonumber
%\\
%&\qquad\times
 e^{\left(\frac{s^{\nu *}}{s_k-s^{\nu *}}\right)}e^{\left(-\frac{s^{\nu *}}{s_0-s^{\nu *}}\right)}.
\label{eq:21*}
\tag{21*}
\end{align}

Remark that  \eqref{eq:21} is different from the correct solution \eqref{eq:21*} to the ODE \eqref{eq:20}. In fact, if one chooses $C\cong 0$ and neglects the second term on the LHS of \eqref{eq:b} while solving ODE (20), one gets \eqref{eq:21}. In Section \ref{s:consequences}, we substantiate that $C$ cannot be simply neglected under the parametric setting of \cite{paper}. Moreover, the removal of the second term leads to an unreliable estimation of the rate of  convergence.

\subsection{Technical Flaw in the Algorithmic Design}\label{ss:limitation}

\begin{figure}[!t]
\centering
\subfloat[$\nu=1/2$]{\includegraphics[width=0.2\textwidth]{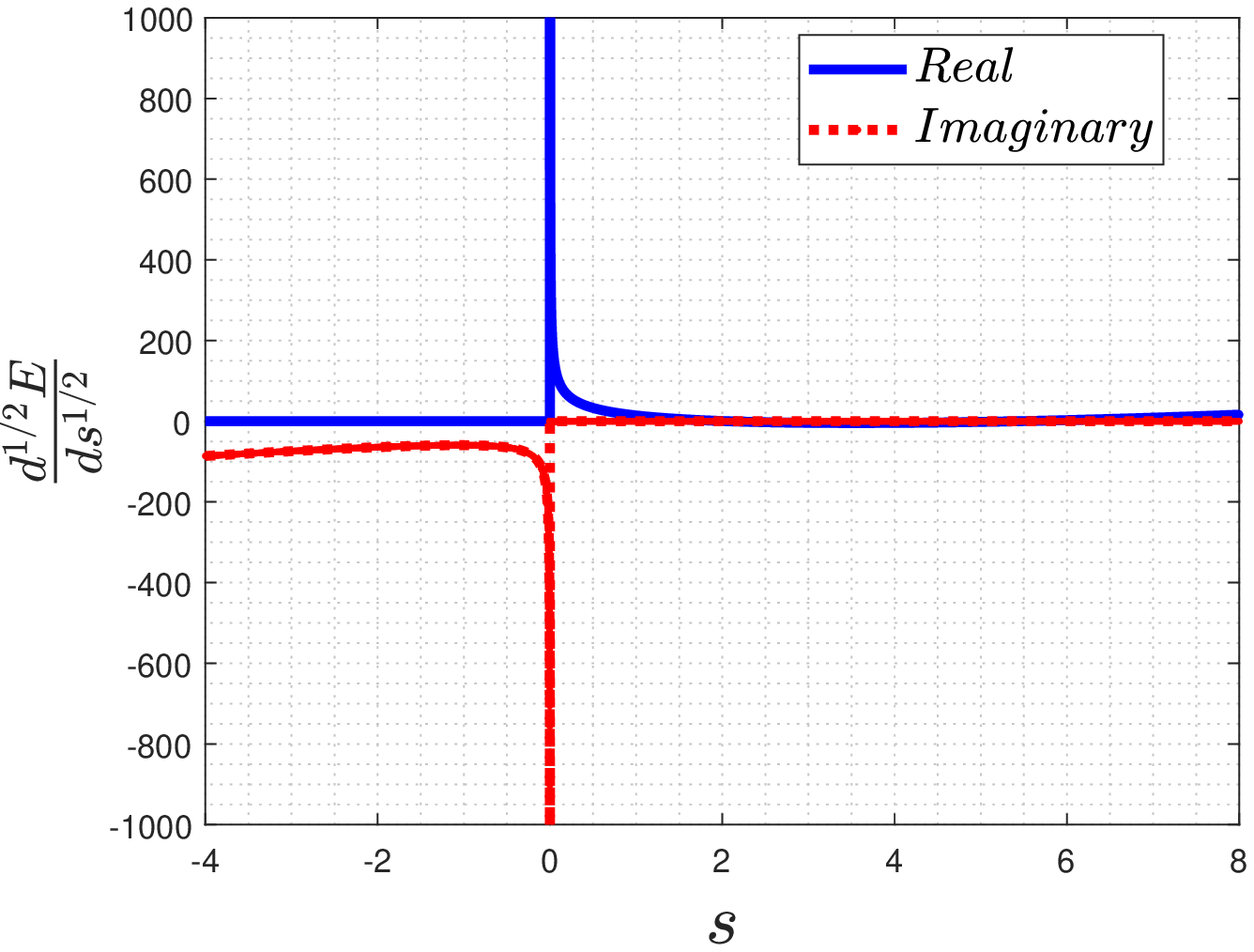}
\label{fig:FD1}}
\hfil
\subfloat[$\nu=3/2$]{\includegraphics[width=0.2\textwidth]{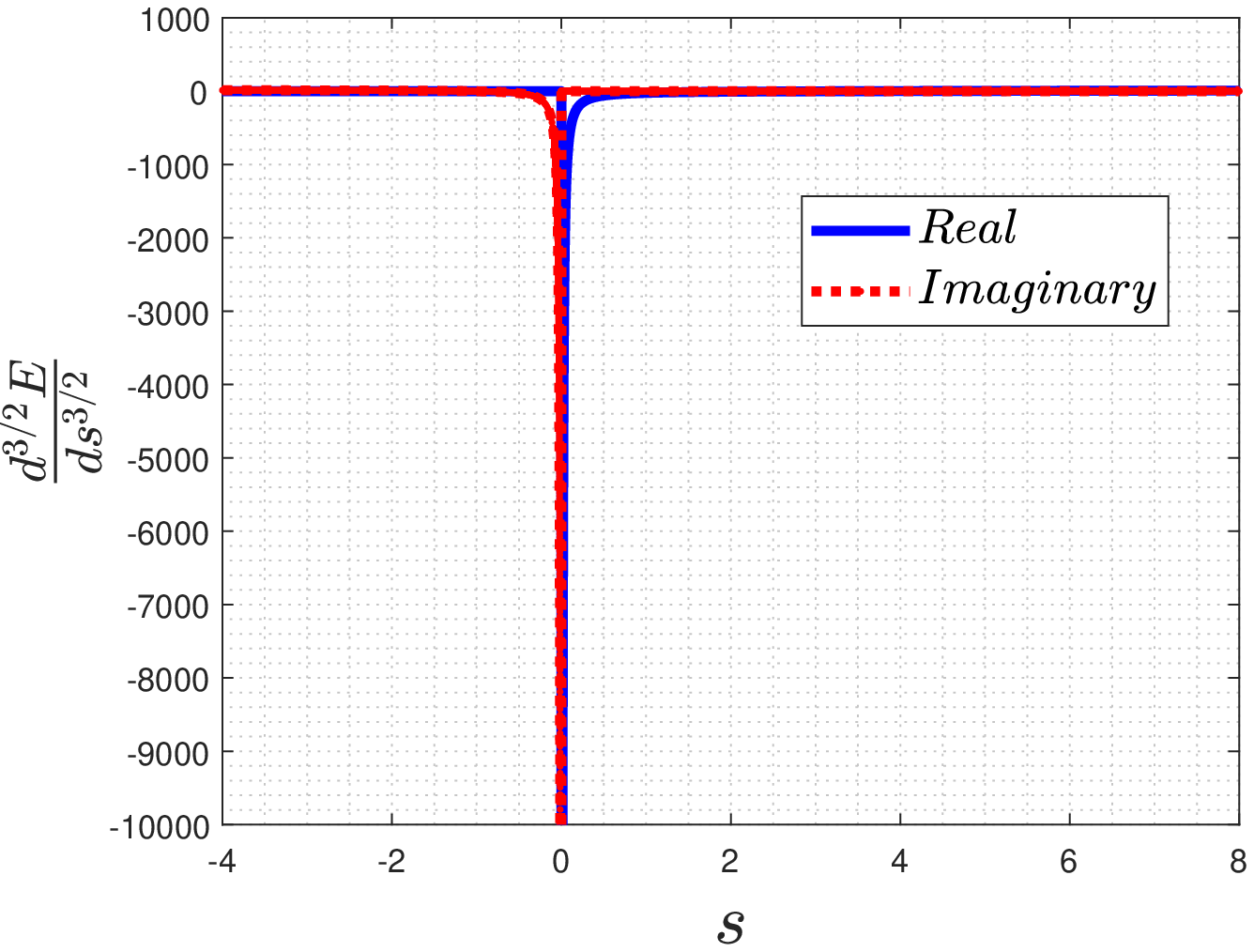}
\label{fig:FD2}}
\caption{Fractional derivative of the energy norm.}
\label{fig:FD}
\end{figure}

The FAL approach in \cite{paper} is proposed for seeking a minimizer $s^{\nu *}$ of the energy norm \cite[Eq. (6)]{paper} in the real domain $\mathbb{R}$.  Both negative and positive minimizers are sought in \cite{paper}. However, the update equation \eqref{eq:19} of the FAL algorithm contains a fractional power of $s_k$ which  becomes complex whenever  $s_k<0$. In particular, ${d^{\nu} E}/{ds^{\nu}}$ for $\nu = 1/2$ and $3/2$ is pure imaginary. In this situation, $s_{k+1}$ will be complex since \eqref{eq:19} is also derived from ${d^{\nu} E}/{ds^{\nu}}$. %Therefore, all the subsequent terms are expected to be complex. 
Consequently, the FAL method is not expected to converge to a real value. In order to elaborate on this point, we evaluate ${d^{\nu} E}/{ds^{\nu}}$ (based on \cite[Eq. (8)]{paper})  using the same parameters as in \cite[Sect. IV-B]{paper}, i.e., we set $E^1_{\min}=10$, $\eta=2$, and $s^{1,*}=5$, $1<\nu\leq 2$, and the domain $-4<s<8$ as used for \cite[Figs. 2(e), 2(d)]{paper}. Then, for  $\nu=3/2$,
\begin{align}
\frac{d^{3/2} E}{ds^{3/2}}
%=&\frac{\left[E^1_{\min}+\eta(s^{1*})^2\right]s^{-{3/2}}}{\Gamma(-{1/2})}-\frac{2\eta s^{1 *} s^{{-1/2}}}{\Gamma({1/2})}+\frac{2\eta s^{{1/2}}}{\Gamma({3/2})}
%\nonumber
%\\
%=&\frac{\left[10+2(5)^2\right]}{-2\sqrt{\pi}}s^{-3/2}-\frac{2(2) (5) }{\sqrt{\pi}}s^{-1/2}+\frac{2(2) }{\sqrt{\pi}/2}s^{1/2}
%\nonumber
%\\
=&-\frac{1}{\sqrt{\pi}}\left(30s^{-3/2}+20s^{-1/2}-8s^{1/2}\right),
\label{eq:7'}
\tag{2'}
\end{align}
which contains fractional powers of $s\in (-4, 8)$.  In particular, at $s=-1$,
\begin{align}
\frac{d^{3/2} E}{ds^{3/2}}\Bigg|_{s=-1}
%=&\frac{-30(-1)^{-3/2}-20(-1)^{-1/2}+8(-1)^{1/2}}{\sqrt{\pi}}
= -\frac{2\iota}{\sqrt{\pi}},
\label{eq:8'}
\tag{3'}
\end{align}
where $\iota = \sqrt{-1}$.  Similarly, the $1/2-$order derivative of the energy norm (based on \cite[Eq. (8)]{paper}) can be calculated as  
\begin{align}
\frac{d^{1/2} E}{ds^{1/2}}
=&\frac{4}{3\sqrt{\pi}}\left(45s^{-1/2}-30s^{1/2}+4s^{3/2}\right),
\label{eq:9'}
\tag{4'}
\end{align}
with parameters as in \cite[Fig. 2(a)]{paper}. Especially, at $s=-1$, 
\begin{align}
\frac{d^{1/2} E}{ds^{1/2}}\Big|_{s=-1}
=&-\frac{316\iota}{3\sqrt{\pi}}.
\label{eq:10'}
\tag{5'}
\end{align}
As a result, \eqref{eq:19} is also complex since it is  based on the same expression of the fractional derivative. Consequently, the future updates $s_{k+1}$ will be complex and the algorithm will not converge to a real value as anticipated. 
%Hence, the statement: ``\emph{If  $s<0$ and $\nu=1/2$, $(d^\nu E/ds^\nu)\equiv 0.$''} (see \cite[Page 658]{paper})  is incorrect in view of \eqref{eq:10'}. 
In order to substantiate this, we plotted the expressions \eqref{eq:d} and \eqref{eq:9'} in Fig. \ref{fig:FD} over the domain $(-4,8)$ using same parameters as in \cite[Fig. 2]{paper}. It is observed that ${d^{\nu} E}/{ds^{\nu}}$ is real as long as $s>0$ and is pure imaginary for $s<0$. Note also that ${d^{\nu} E}/{ds^{\nu}}$ is singular at $s=0$ which actually justifies  that $s_0=0$. 

\section{Discussion}\label{s:consequences} 

\subsection{Reliability of the Rate of Convergence}\label{ss:reliability}
Let us discuss some consequences of the flaws indicated in Section \ref{ss:error}. First, it is worthwhile precising that the FAL approach is based on left Riemann-Liouville fractional derivative \cite[Eq. (3)]{paper} (instead of Gr\"{u}nwald-Letnikov derivative as pretended in \cite{paper}) with $a=0$ . Therefore, FAL is valid only for $s>0$ and $s_0=0$. Consequently, Eq. \eqref{eq:d} suggests that $C\cong \ln |s^{\nu*}|+1$. Since $s^{\nu*}$ is unknown sought value, one cannot simply set $C\cong 0$ in \eqref{eq:b} to get \eqref{eq:21}.  

On the other hand, the approximation \eqref{eq:21}, derived from \eqref{eq:21*} by ignoring $\exp\left({s^{\nu *}}/{(s_k-s^{\nu *})}\right)$ and choosing $C\cong 0$, is highly unreliable. The convergence analysis in \cite{paper} is based entirely on the estimate \eqref{eq:21}. By choosing $\mu$ such that 
\begin{align}
\lim_{k\to+\infty} k \chi=+\infty,\quad\text{with}\quad \chi:= \frac{2\mu\eta}{\Gamma(3-\nu)\nu \left(s^{*\nu}\right)^{\nu-1}}>0,
\label{eq:f}
\tag{6'}
\end{align}
it is suggested in \cite{paper} that the \emph{algorithm converges at the rate $\exp(-\chi k )$}. In fact,  since $(s_k)_{k\in\mathbb{N}}$ is assumed to be convergent to $s^{\nu*}$,   $(s_k-s^{\nu*})\to 0$ as $k\to +\infty$. Hence, $s^{\nu *}/(s_k-s^{\nu *})\to+\infty$ and consequently, $\exp\left(s^{\nu *}/(s_k-s^{\nu*})\right)\to +\infty$ when $s^{\nu *}/(s_k-s^{\nu *})$ is positive and $k\to+\infty$.  
Therefore, the  product $\exp\left(-\chi k\right) \exp\left(s^{\nu *}(s_k-s^{\nu *})^{-1}\right)$  has an indeterminate form $0\times \infty$. One cannot guarantee that it will approach to $0$. 
%$\exp\left(-\chi k\right) \exp\left(s^{\nu *}(s_k-s^{\nu *})^{-1}\right)\to 0$. 
Even if it does so, the factor $\exp\left(s^{\nu *}(s_k-s^{\nu *})^{-1}\right)$ will severely impede the decay of $\exp\left(-\chi k\right)$, which will be grandiloquent as the rate of convergence of FAL. 

\begin{figure}[!t]
\centering
\subfloat[$\chi=0.25$]{\includegraphics[width=0.2\textwidth]{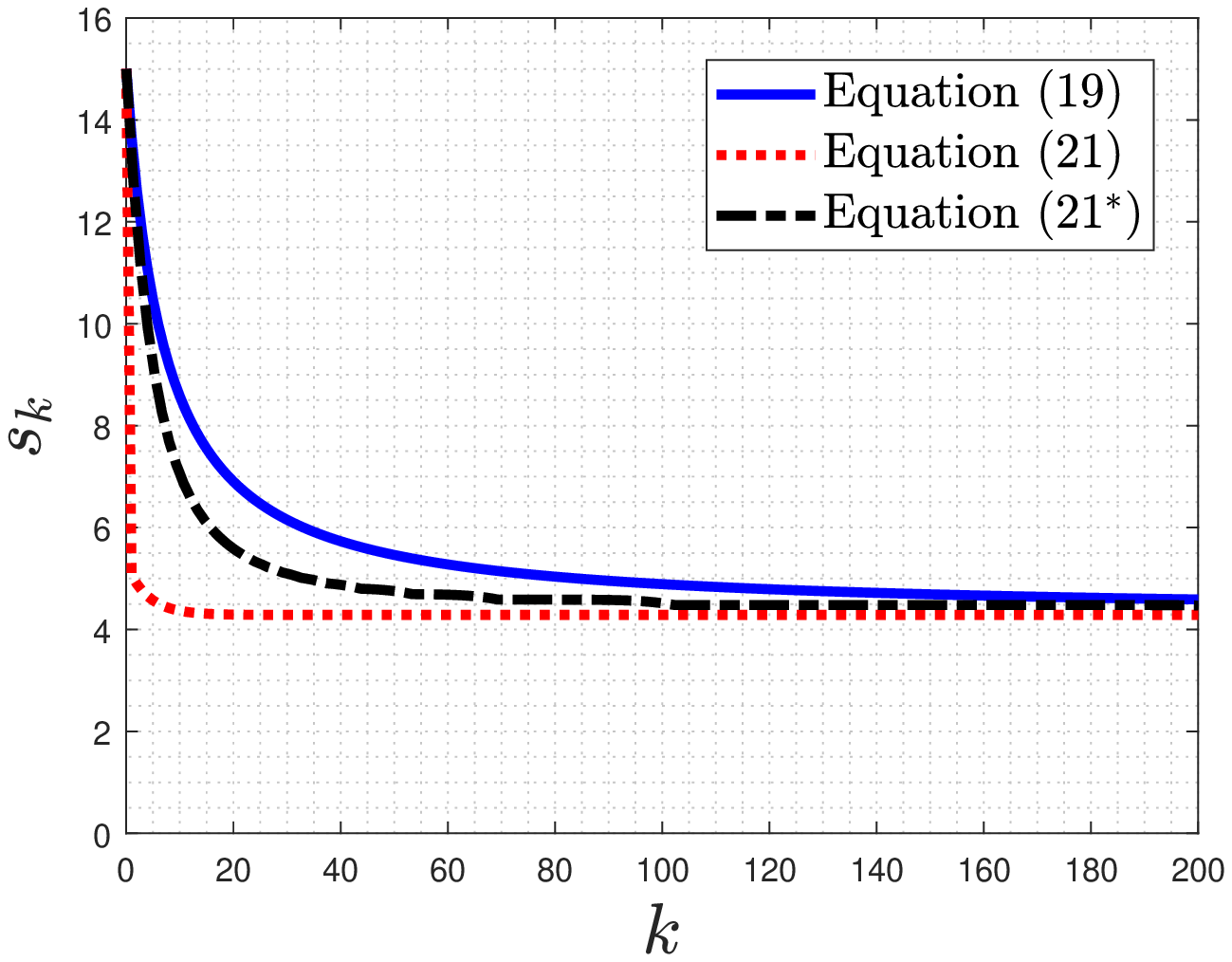}
\label{fig:rates1}}
\hfil
\subfloat[$\chi=1.75$]{\includegraphics[width=0.2\textwidth]{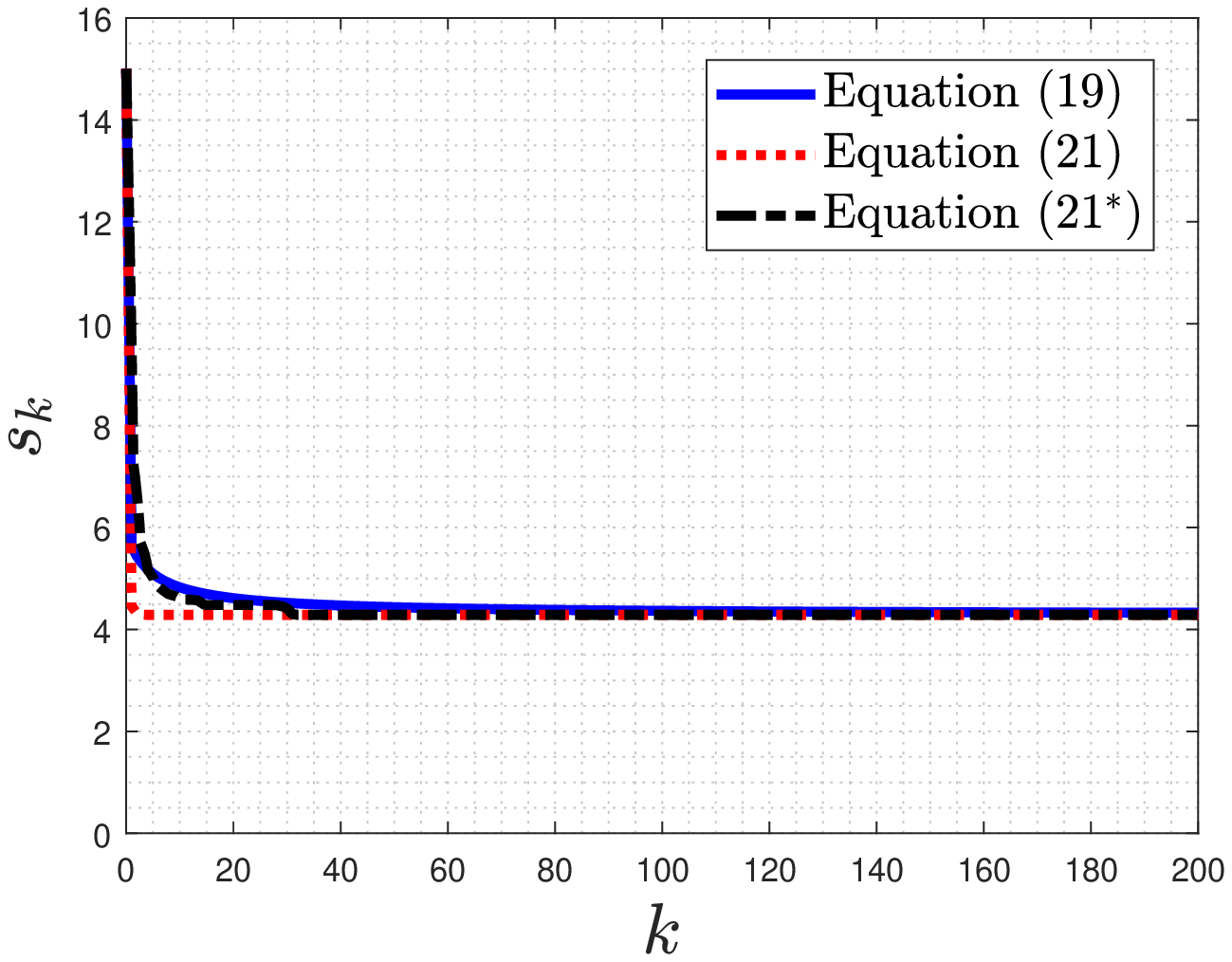}
\label{fig:rates2}}
\caption{Estimation of the rate of convergence.}
\label{fig:rates}
\end{figure}

In order to elaborate on this point, we have compared the rates of convergence based on estimates \eqref{eq:19}, \eqref{eq:21}, and \eqref{eq:21*} in Fig. \ref{fig:rates}. We choose same parameters as in \cite[Fig.5(a)]{paper}. The computational results indicate that the FAL (with update rule \eqref{eq:19}) converges at a very slow rate as compared to that predicted by \eqref{eq:21}. When $\chi=0.25$,  \eqref{eq:21} suggest that FAL converges to the sought value $s^{\nu *}=4.2856$ after only 29 iterations with $s_{29}\approx 4.2856$. On contrary,  \eqref{eq:19} suggests that after $k=1948$ iterations $s_k\approx 4.316$. On the other hand, \eqref{eq:21*} predicts that a steady state is achieved at $k=414$ with $s_{414}\approx 4.2856$. 
Similarly, when $\chi=1.75$, the actual number of iterations for FAL  to achieve a steady state is $k=1741$ whereas \eqref{eq:21} and \eqref{eq:21*} predict $k=5$ and $k=56$, respectively.

 Following remarks are in order. First, \eqref{eq:21} does not provide any reliable estimate for the rate of convergence as the actual convergence is roughly two orders of magnitude slower than the predicted rate. Second, \eqref{eq:21*} also predicts a convergence  almost an order of magnitude faster than the actual rate, yet, it provides much superior estimation than \eqref{eq:21}. Thirdly, based on these observations, it seems inappropriate to consider $s_k$ as a discrete sample of a continuous function $s(t)$ based on which both \eqref{eq:21} and \eqref{eq:21*} are derived.  
%\begin{figure}[!t]
%\centering
%\includegraphics[width=0.23\textwidth]{CASE1}
%\quad
%\includegraphics[width=0.23\textwidth]{CASE2}
%\caption{Simulation results for the network.}
%\label{fig:rates}
%\end{figure}

\subsection{Consequences of the Flaw in Algorithmic Design}\label{ss:conseq}

In view of the remarks in  Section \ref{ss:limitation}, it is clear that for negative sought values, the FAL update weight $s_{k+1}$ in \eqref{eq:19} becomes complex and cannot converge to a real negative desired output.  As mentioned above, the fractional gradient (\cite[Eq. (8)]{paper})   is valid over the domain $(0,s)$. Therefore, the algorithm cannot be used for negative values of the independent variable. This is the main reason that the fractional derivative appears to be complex for $s<0$. As a consequence, almost every simulation in \cite{paper} is affected and is unreliable. 
\begin{enumerate}
\item \cite[Figs. 2(a), (b), (d), (e), (g), and (h)]{paper} are counter-factual as the fractional derivative in all these cases is complex. Particularly for $\nu=1/2$ or  $3/2$, the fractional derivative is pure imaginary (see, for example, \eqref{eq:d}-\eqref{eq:10'} or Fig. \ref{fig:FD} in this note). 
%Note that in all these examples, when $\nu=1/2$ or $\nu=3/2$, the derivatives of the energy is identically zero over all $s<0$ (refer to\cite[Fig. 2]{paper} on contrary to Fig. \ref{fig:FD}). It is suspected that only the real part of the derivative is plotted. 
For $\nu=0$  (i.e., no derivative is taken), the quadratic energy function \cite[Eq. (6)]{paper} is expected to have a parabolic graph. However, in \cite[Fig 2(h)]{paper}, it appears to be a straight line, which is impossible. Similar  observations also hold for \cite[Figs. 2(b), (d), (e), and (g)]{paper}. 

\item In \cite[Fig. 3(b)]{paper}, the fractional derivative is evaluated over the domain $s<0$. Therefore,  derivative should be complex valued. 

\item In \cite[Fig. 4(b)]{paper}, a negative optimal value  $s_2^{\nu *}=-0.6406$ is sought. In fact, with initial step $s_{20}=-0.25$ and the  parameters for \cite[Fig. 4(b)]{paper}, even $s_{21}$ becomes complex if \cite[Eq. \eqref{eq:19}]{paper} is used. If \cite[Eq. \eqref{eq:21}]{paper} is used  then the exponent on the RHS becomes complex (due to the  term $(s^{\nu *})^{\nu-1}$). 

\item In \cite[Fig. 5]{paper}, the rate of convergence is evaluated for different choices of  $\mu$ and  $\chi$.  As discussed in Section \ref{ss:reliability}, the displayed results are misleading and grandiloquent (see Fig. \ref{fig:rates} in this note).  Note that $s_0=15$ is assumed for \cite[Fig. 5]{paper} whereas $s_0=0$ is tacitly assumed in the derivation of the FAL algorithm. 

\item The results in \cite[Fig. 6]{paper} are also affected by the complex outputs when the $x$ or the $y$ component is varying over a part of the negative axis as multi-dimensional FAL is essentially a generalization of the 1-D FAL.  
\end{enumerate}

\subsection{Comparison to \cite{Bershad}}
 
In \cite{Bershad}, Bershad, Wen, and Cheung So, have already debated the unsuitability of fractional learning frameworks for adaptive signal processing \cite{Raja}. 
Theoretical obeservations in this note can be compared to those made in \cite{Bershad} through a variety of experimental results (see \cite[Sect. 1 and Remark 1]{Bershad}).  In fact, it is well-known that the LMS algorithm is a stochastic version of the steepest descent algorithm when the statistics of the input are unknown. Thus, \cite[Eq. (1)]{Bershad} can be compared directly to \cite[Eq. (19)]{paper}. 

Based on extensive experiments, the following conclusions have been drawn in \cite[Page 225] {Bershad}.  
\begin{enumerate}
\item The fractional variants of the LMS are only useful when all the update weights are positive but their performance is comparable to that of the LMS. That is, under no conditions fractional variants of LMS perform better than the  standard LMS. 

\item In case when some of the  update weights are negative, the fractional variants of LMS render complex outputs (see \cite[Remark 1]{Bershad}). Moreover, even when the absolute operator is employed in the fractional algorithms (see, for instance, Refs. 3 and 5 in \cite{Bershad}),  their performance is inferior than standard LMS. Finally, if only the real part of the complex update weight is employed, the fractional LMS reduces to LMS  with a slower convergence rate.
\end{enumerate}

Observe that the FAL method proposed in \cite{paper} has similar drawbacks as highlighted in \cite{Bershad} for fractional frameworks for adaptive signal processing. Precisely, as debated in Sections \ref{ss:reliability} and \ref{ss:conseq}, the FAL method has limited applicability for broad spectrum of adaptive learning problems due to complex outputs and has slow convergence rate when the update iterates remain real.

\section{Conclusion}\label{s:conclusion}

In this comment, some serious concerns over the derivation of the rate of convergence of Fractional Adaptive Learning (FAL) approach proposed in \cite{paper} were raised. It is established that the convergence analysis perfomed in \cite{paper} is unreliable in general and the FAL algorithm converges much slower than anticipated. It was also highlighted that the FAL method can practically work only for positive domains. Over negative domains or whenever its iterative update becomes negative, the FAL algorithm furnishes a  complex output  due to the presence of fractional powers in its update rule. In this situation, the algorithm is not expected to converge to a real sought value.  Moreover, thanks to the analogy of the FAL algorithm with fractional variants of Least Mean Squares (LMS) for adaptive signal processing \cite{Raja}, the analysis performed by Bershad, Wen, and Cheung So \cite{Bershad} suggests that FAL is not better than LMS under any condition. Their performances are nearly the same but the FAL approach is much more complicated than LMS. Finally, it is needless to say that the multi-dimensional variant of the FAL also inherits the same flaws and is unreliable.

%\appendices

%\section{}

%\section*{Acknowledgment}


\ifCLASSOPTIONcaptionsoff
  \newpage
\fi

\end{document}